\documentclass[11pt]{elsarticle}
\biboptions{numbers,sort&compress}

\usepackage{amsthm}
\usepackage{geometry}
\usepackage{amsmath,amssymb,amsfonts}
\usepackage{algorithmic}
\usepackage{graphicx}
\usepackage{textcomp}

\newtheorem{theorem}{Theorem}

\usepackage{appendix}
\usepackage{setspace}
\usepackage{multirow}
\usepackage{xcolor}
\usepackage{bbding}
\usepackage{booktabs}
\usepackage{algorithm}
\usepackage{subfigure}
\usepackage{threeparttable}
\usepackage{hyperref}

\journal{POEMA}

\begin{document}

\begin{frontmatter}



\title{Safe Controlled Invariance for Linear Systems Using \\Sum-of-Squares Programming}

\author{Han Wang}\ead{han.wang@eng.ox.ac.uk}
\author{Kostas Margellos}\ead{kostas.margellos@eng.ox.ac.uk}
\author{Antonis Papachristodoulou}\ead{antonis@eng.ox.ac.uk}
\address{Department of Engineering Science, University of Oxford, Oxford, OX1 3PJ, United Kingdom}

\begin{abstract}
Safety is closely related to set invariance for dynamical systems. However, synthesizing a safe invariant set and at the same time synthesizing the associated safe controller still remains challenging. In this note we introduce a simple invariance-based method for linear systems with safety guarantee. The proposed method uses sum-of-squares programming.

\end{abstract}




\end{frontmatter}


\section{Introduction} \label{sec:intro}
In this note we show how polynomial optimisation can be used to analyze safety for linear systems. A dynamical system is said to be safe with respect to a \emph{safe set} $\mathcal{S}$ if all the trajectories starting from an \emph{initial set} $\mathcal{I}\subseteq \mathcal{S}$ cannot leave $\mathcal{S}$. Verifying this property requires the solution of a reachability optimal control problem \cite{margellos2011hamilton}, where a Hamilton-Jacobi-Issac partial differential equation is formulated. An alternative verification method is known as the method of \emph{Barrier Certificates}, which checks a condition on the boundary of the safe set under system dynamics \cite{prajna2004safety}. A barrier certificate is a function whose the zero-sub (or zero-super in some recent works \cite{ames2016control} \cite{wang2022safety}) level set is invariant. The boundary condition determines the sub-tangentiality of the set, thereby is a sufficient condition to invariance according to Nagumo's Theorem \cite{blanchini1999set}. For polynomial systems with semi-algebraic safe/initial sets, the conditions for the existence of a barrier certificate can be formulated as a sequence of constraints, which can be convexified and solved using \emph{sum-of-squares} programming \cite{parrilo2000structured}. The resulting problem is a semi-definite program, and solvers such as SeDuMi \cite{doi:10.1080/10556789908805766} employing interior-point methods can provide a numerical solution efficiently.

Although sum-of-squares programs and barrier certificates are shown to be reliable and efficient in verifying safety for general polynomial systems, there are two problems encountered.
\begin{itemize}
    \item For controlled dynamical systems, the synthesis optimisation problem is bilinear due to the presence of the additional control input.
    \item The sum-of-squares programs do not scale well with the size of the dynamical systems considered.
\end{itemize}
Solving these two problems for general polynomial systems is difficult and beyond the scope of this note. Motivated by a sum-of-squares method to construct Lyapunov functions \cite{papachristodoulou2002construction} and linear system set invariance \cite{blanchini2008set}, we focus on safe invariance set construction and feedback controller design for linear systems. A relatively simple and scalable method based on ellipsoidal invariant sets is proposed.

\section{Main Results}\label{sec:conclusion}
We consider a continuous-time linear time-invariant system
\begin{equation}\label{eq:linsys}
    \dot x=Ax+Bu,
\end{equation}
where $x(t) \in\mathbb{R}^n$ denotes the state, and $u(t) \in\mathbb{R}^m$ the control input, while the matrices $A\in\mathbb{R}^{n\times n}$ and $B\in\mathbb{R}^{n\times m}$. Define the compact semi-algebraic set $\mathcal{S}:=\{x|s(x)\ge 0\}$ as the safe set. Throughout the paper, we assume that the system can be safely controlled within $\mathcal{S}$. For the sake of brevity, we assume the system is fully observable without state and output noise. Also, it is natural to assume that $(A,B)$ is controllable. The control input $u$ is bounded by $||u||_2^2\le u_{\max}$, where $u_{\max}>0$ is a positive scalar.

The system is assumed to be locally stabilizable around $x^*\in\mathcal{S}$. Among different types of invariance sets, we select the ellipsoidal invariance set as a candidate. Our choice stems from the fact that, if the system \eqref{eq:linsys} is stabilizable (or locally stabilizable) with state feedback control law $u=kx$, there exist a quadratic Lyapunov functions $V(x)=x^\top Px$ \cite{khalil2002nonlinear}, Theorem 4.17. Then every positive-sub level set of $V(x)$ serves as an invariant set. Under this set-up, we use 
\begin{equation}\label{eq:barrier}
    b(x)=-x^\top Px+l,
\end{equation}
where $P$ and $l$ are parameters to be determined for a candidate function for the safe invariance set.

According to the Nagumo's theorem and barrier certificates conditions \cite{prajna2004safety}, $b(x)$ satisfies:
\begin{subequations}\label{eq:condition}
\begin{align}
    &\mathrm{for~any}~x~\mathrm{such}~\mathrm{that}~b(x)=0,\exists u\in\{u|||u||_2^2\le u_{\max}\},~\frac{\partial b(x)}{\partial x}(Ax+Bu)\ge 0,\label{eq:condition1}\\
    &b(x)<0~\mathrm{for~any}~x~\mathrm{such}~\mathrm{that}~s(x)<0.\label{eq:condition2}
\end{align}
\end{subequations}
To guarantee the condition \eqref{eq:condition1} holds, we have
\begin{equation}\label{eq:sup}
    \sup_{||u||_2^2\le u_{\max}} \frac{\partial b(x)}{\partial x}(Ax+Bu)\ge 0, \mathrm{for~any}~x~\mathrm{such}~\mathrm{that}~b(x)=0.
\end{equation}
By \eqref{eq:barrier}, we have $\frac{\partial b(x)}{\partial x}=2Px$. The maximizing $u$ satisfies
\begin{equation}\label{eq:inner}
    <B^\top\frac{\partial b(x)}{\partial x}^\top,u>=||\frac{\partial b(x)}{\partial x}B||\cdot ||u||,
\end{equation} 
where $<\cdot,\cdot>$ denotes the inner product operator. We directly have $u=-\zeta B^\top P^\top x$, $\zeta$ is a positive scaled scalar. Without control limitation, $\zeta$ should be large enough to maximize $\frac{\partial b(x)}{\partial x}Bu$. Under the assumption that $u$ is bounded, $\zeta$ is a function of both state $x$ and matrix $P$. Condition \eqref{eq:condition1} turns out to be
\begin{equation}\label{eq:newcondition1}
    x^\top(-PA-A^TP+2\zeta PBB^\top P^\top)x\ge 0,
\end{equation}
which is equivalent to 
\begin{equation}\label{eq:newcondition2}
    -PA-A^TP+\zeta PBB^\top P^\top \succeq 0.
\end{equation}
We omit the boundary condition $b(x)=0$ since $\mathcal{B}$ is convex and compact \cite{blanchini2008set}. In fact, for a point $y$ such that $b(y)\ne 0$, we can always find $x\in\mathbb{R}^n$ and $t\in\mathbb{R}$, such that $x=ty$ and $b(x)=0$. We directly have that $x^\top(-PA+\zeta PBB^\top P^\top)x\ge 0 \Leftrightarrow y^\top(-PA+\zeta PBB^\top P^\top)y\ge 0. $

To overcome the bilinear term $\zeta PBB^\top P^\top$ in \eqref{eq:newcondition2}, we propose the following relaxed conditions:
\begin{subequations}\label{eq:relaxedconditions}
\begin{align}
&-PA-A^TP+2\zeta \hat P\succeq 0,\label{eq:rc1}\\
& \left[ {\begin{array}{*{20}{c}}
{\hat P}&{PB}\\
{{B^ \top }{P^ \top }}&I
\end{array}} \right]\succeq 0.\label{eq:rc2}
\end{align}
\end{subequations}
Theorem \ref{th:suff} states that \eqref{eq:relaxedconditions} is sufficient for \eqref{eq:newcondition2}.
\begin{theorem}\label{th:suff}
If there exist $\zeta>0$, $P\succeq 0$ and $\hat P$ satisfying \eqref{eq:relaxedconditions} then $P$ and $\zeta$ also satisfy \eqref{eq:newcondition2}.
\end{theorem}
If we fix $\zeta$, then the new relaxed conditions \eqref{eq:rc1} and \eqref{eq:rc2} are both convex semi-definite constraints. Now we give the overall synthesis approach using the following two sum-of-squares programs.
\begin{subequations}\label{eq:sos1}
\begin{align}
    \mathop{\max_{P,\hat P,l,\sigma}}&~~~l~~~~~~~~~~~~~~~~~~~~~~~~~~~~~\\
    \mathrm{such~that}&~~~x^\top Px-l+\sigma s(x)\in\Sigma[x],~\eqref{eq:relaxedconditions} ~~~~~~~~~~~~~~~~~~~~~~~~~~~~~
\end{align}
\end{subequations}
\vspace{-7ex}
\begin{subequations}\label{eq:sos2}
\begin{align}
    \mathop{\max_{\zeta^2,\sigma}}&~~~\zeta^2\\
    \mathrm{such~that}&~~~    -\zeta^2 x^\top PBB^\top P^\top x-\sigma (-x^\top Px+l)+ u_{\max}\in\Sigma[x],
\end{align}
\end{subequations}
Here we use $\Sigma[x]$ to represent the set of sum-of-squares polynomials in $x$, and $\sigma\in\Sigma[x]$ to be a sum-of-squares multiplier. Details of sum-of-squares program are referred to \cite{papachristodoulou2005tutorial}. The programs \eqref{eq:sos1} and \eqref{eq:sos2} are computed iteratively. Solving \eqref{eq:sos1} requires to fix $\zeta$ for convexity, hereby $\zeta$ is initialized to be a small constant. We remark that the program \eqref{eq:sos2} uses $\zeta^2$ as the variable rather than $\zeta$ to avoid the trivial bilinearity of $\zeta\cdot\zeta$.
 
\begin{figure*}[h]
    \centering
        \subfigure[Phase portrait for the system]{
        \includegraphics[width=0.47\textwidth]{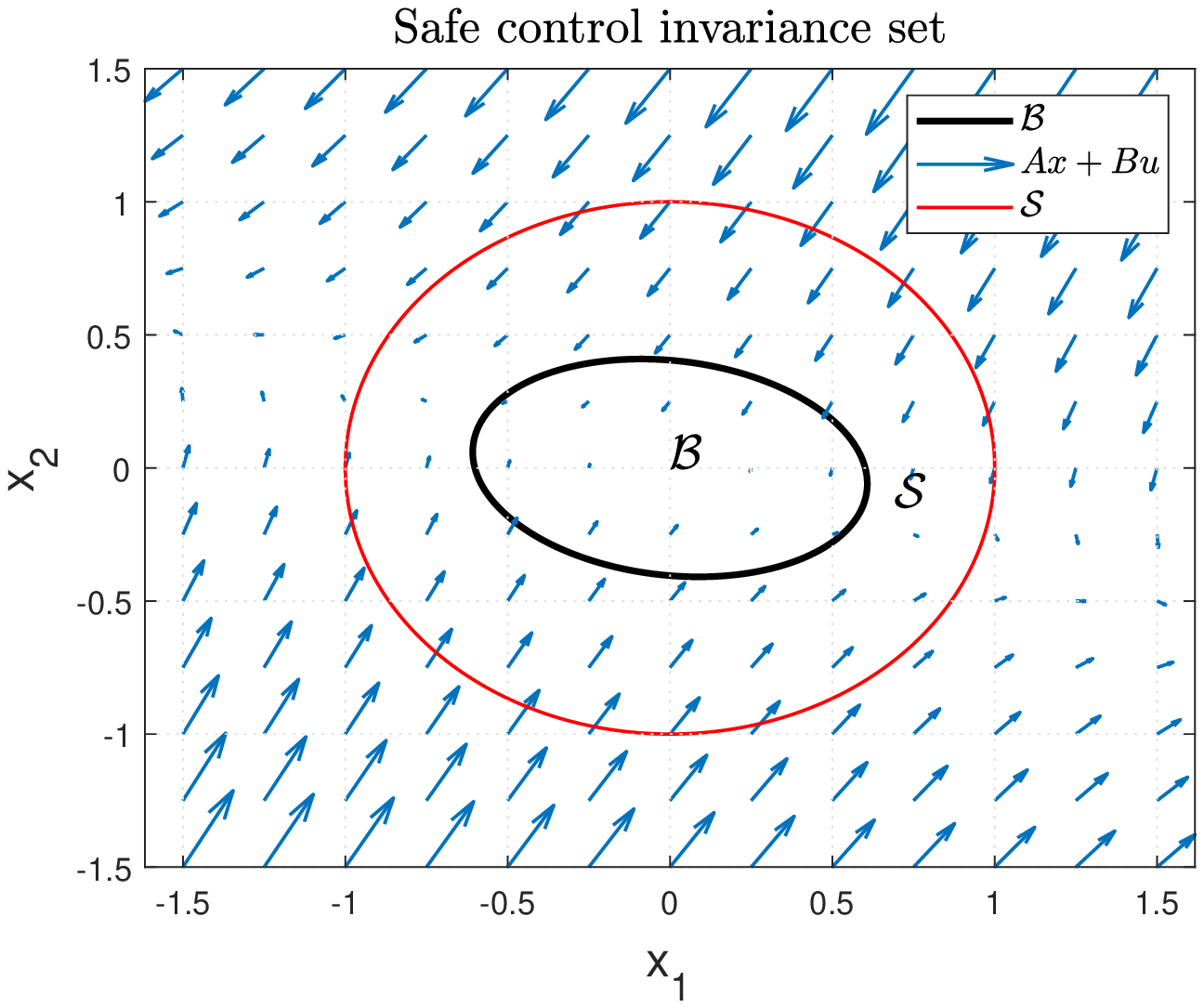}
        \label{fig:Nonlinsys}
    }
    \subfigure[Level set of $||u||_2^2$, $u=-\zeta B^\top P^\top x$]{
        \includegraphics[width=0.47\textwidth]{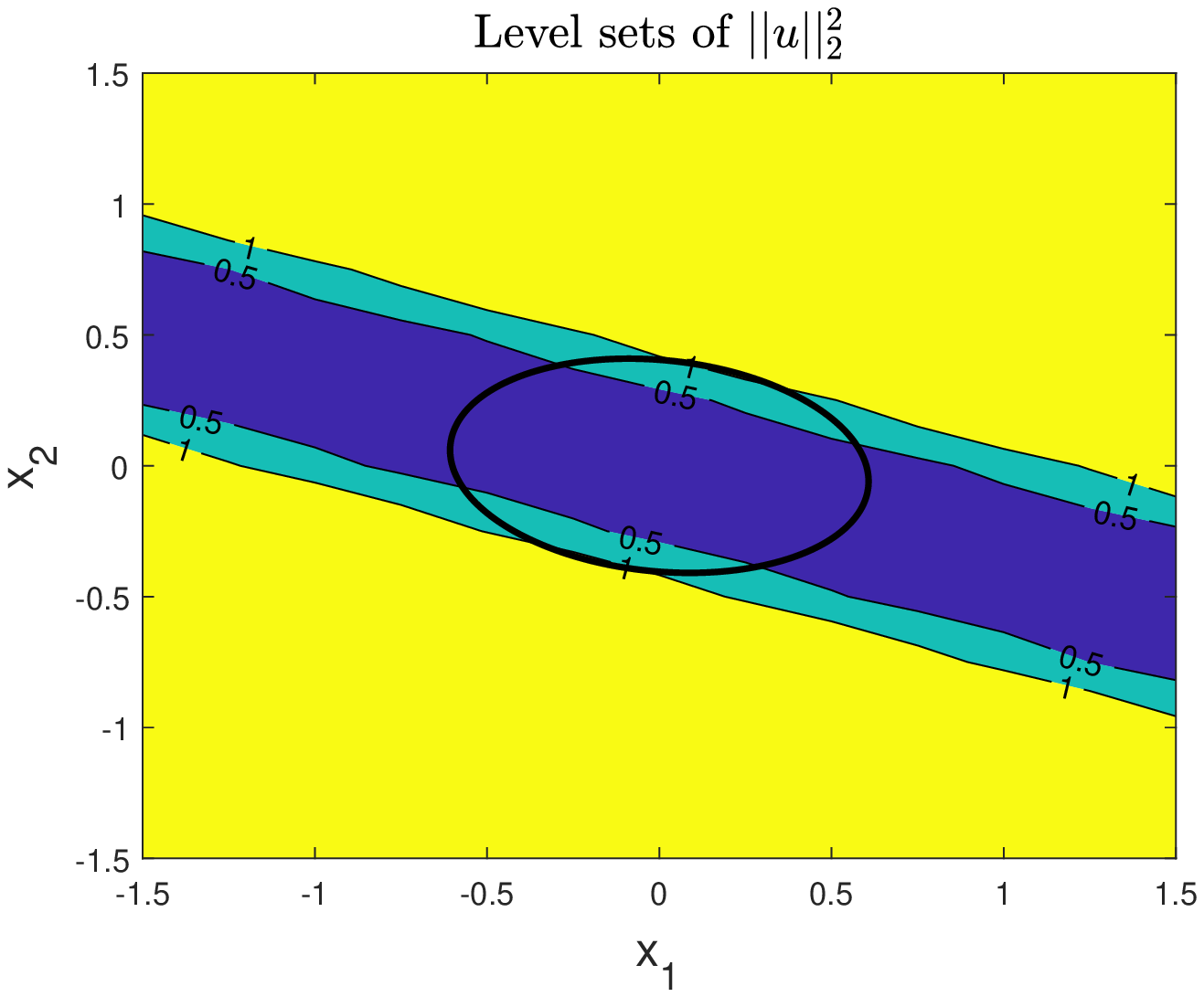}
        \label{fig:NonSolU1}
    }
    \caption{The interior of the red disc represents the safe set, the black closed curve inside $\mathcal{S}$ is the control invariance set $\mathcal{B}$. The arrows in the figure represent the vector field.} 
    \label{fig:NonLinSysCon}
\end{figure*}

We then test our results on a second order linear system with $A=\left[ {\begin{array}{*{20}{c}}
{0.8}&{0.7}\\
{ - 0.4}&{ - 0.6}
\end{array}} \right]$, $B=\left[ {\begin{array}{*{20}{c}}
1&1\\
1&1
\end{array}} \right]$, $u_{\max}=1$. The safe set is the interior of a disc: $\mathcal{S}:=\{x|-x_1^2-x_2^2+1\ge 0\}$. Without control input, the system is an unstable spiral which certainly breaks safety requirement inside the closed safe set. Figure \ref{fig:NonLinSysCon} shows the synthesized invariant set and value of control input. The sum-of-squares programs \eqref{eq:sos1} and \eqref{eq:sos2} are solved using SOSTOOLS v400 \cite{sostools} with SeDuMi v1.3.5 \cite{doi:10.1080/10556789908805766} as the SDP solver.

\section{Conclusion and Future Directions}
In this note we propose a simple method to synthesize the safe controlled invariant set with sum-of-squares programming. The method uses semi-definite relaxation to overcome part of the bilinearity introduced by the presence of control input. In the future work we aim at accounting for the presence of control input constraints and encoding the overall problem as a convex program.

\bibliographystyle{elsarticle-num} 
\bibliography{ref.bib} 




\end{document}